\documentclass[a4paper,12pt]{article}
   \usepackage{hyperref}
   \usepackage{amsmath}
   \usepackage{amsfonts}
   \usepackage{amssymb}
   \usepackage{mathrsfs}
\newcommand{\R}{\mathbb{R}}

\newtheorem{thm}{Theora}[section]
\newtheorem{Theo}[thm]{Theorem}

\newtheorem{Lem}[thm]{Lemma}

\newtheorem{Rem}[thm]{Remark}

\newcommand{\be}{\begin{equation}}
\newcommand{\ee}{\end{equation}}
\newcommand{\ba}{\begin{array}}
\newcommand{\ea}{\end{array}}
\newcommand{\bea}{\begin{eqnarray*}}
\newcommand{\eea}{\end{eqnarray*}}
\newcommand{\bean}{\begin{eqnarray}}
\newcommand{\eean}{\end{eqnarray}}

\makeatletter \@addtoreset{equation}{section}

\makeatother
\begin{document}
\title{\bf{The Ratio of Eigenvalues of the Dirichlet Eigenvalue Problem for Equations with One-Dimensional
p-Laplacian}}

\author{Jamel Ben Amara  and  Jihed Hedhly}
\date{}
\date{}
\date{}
\maketitle{ABSTRACT}. Chao-Zhong Chen et al. $[{Proc}.$ ${Amer.
Math. Soc},2013],$ proved the upper estimate
 $\frac{\lambda _{n}}{\lambda _{m}}\leq \frac{%
n^{p}}{m^{p}}$ $ (n>m\geq 1) $ for Dirichlet Shr\"{o}dinger
operators with nonnegative and single-well potentials. In this paper
we discuss the case of nonpositive potentials $q(x)$ continuous on
the interval $[ 0,1] $. We prove
that if $q(x)\leq 0$ and single-barrier then $\frac{\lambda _{n}}{\lambda _{m}}\geq \frac{n^{p}%
}{m^{p}}$ for $\lambda _{n}>\lambda _{m}\geq -2q^{\ast },$ where
$q^{\ast}=\inf\{q(0), q(1)\}$. Furthermore, we show that there
exists $\ell_{0}\in ( 0,1] $ such that for all
$\ell\in(0,\ell_{0}],$ the associated eigenvalues $(\lambda
_{n}(\ell)) _{n\geq 1}$ (of the problem defined on $[0,\ell]$)
satisfy  $ \lambda _{1}( \ell)>0$ and $\frac{\lambda _{n}(
\ell)}{\lambda _{m}( \ell) }\geq \frac{n^{p}}{m^{p}}$ $n>m\geq 1$.
The value $\ell _{0}$ satisfies the following estimate
$0<\ell_{0}\leq \sqrt[p]{\frac{-p}{3q^{*}}}$.

~~~~~~~~~~\\
2010 Mathematics Subject Classification. Primary 34A34, 34L15.\\
Key words and phrases. p-Laplacian, eigenvalue ratio,
single-barrier, generalized Pr\"{u}fer substitution.

\section{Introduction}

We consider the eigenvalue problem for the one-dimensional
p-Laplacian
\begin{equation}\label{1.1}
-((y^{\prime }) ^{( p-1)})^{\prime }=(p-1)(\lambda
-q(x))y^{(p-1)},x\in[ 0,1],
\end{equation}
with $f^{( p-1) }=\vert f\vert ^{p-1}sgn(f)=\vert f\vert ^{p-2}f,~
p>1$ and $q$ nonpositve continuous potential in $[ 0,1]$, with
Dirichlet boundary conditions
\begin{equation} \label{1.2}
y(0)=y(1)=0
\end{equation}
It is known $({see} ~~ \cite{14,15}) $ that the spectrum of problems
$\eqref{1.1} -\eqref{1.2}$ consists of a growing sequence of
infinitely point $\lambda _{1}<\lambda _{2}<....<\lambda _{n}...$.\\
Recently, the issues of optimal estimates for the eigenvalue ratios $\frac{\lambda _{n}%
}{\lambda _{m}}$  (where $ n>m$) for the problem
$\eqref{1.1}-\eqref{1.2}$ in the case $p=2$ have attracted a lot of
attention $($cf. $\cite{1}$, $\cite{HL}$, $\cite{4}$, $\cite{5}$$)$.
In particular, M. Horv\'{a}th and M. Kiss $\cite{5}$ showed that for
nonnegative single-well potentials $\frac{\lambda _{n}}{\lambda
_{m}}\leq \frac{n^{2}}{m^{2}}.$ Let $0\leq x_{0}\leq 1$ be fixed.
Following Ashbaugh and Benguria $\cite{2}$ we call the function $q$
a single-well function if $q$ is decreasing in $[0,x_{0}]$ and
increasing in $[x_{0},1].$ Analogously, $q$ is called a
single-barrier function if it is increasing in $[0,x_{0}]$ and
decreasing in $[x_{0},1].$ Recently the authors  $\cite{6}$ showed,
among other results, that for nonpositive single barrier potentials
\begin{equation}\label{JH}
\frac{\lambda _{n}}{\lambda _{m}}\geq \frac{n^{2}}{m^{2}},\ \ \
\lambda _{n}>\lambda _{m}\geq-2q^{\ast}.
\end{equation}
In the case $p>1$, Gabriella Bogn\'{a}r and Ond\~{r}ej Do\~{s}ly
$\cite{13}$ proved an optimal upper estimate $\frac{\lambda
_{n}}{\lambda _{m}} \leq \frac{n^{p}}{m^{p}}$ of Dirichlet
Schrodinger operator with nonnegative differentiable and single-well
potentials $q(x)$. Later, Chao-Zhong Chen et al. $\cite{12} $
extended the results in $\cite{13}$ without assumption conditions on
the differentiability of $q(x)$. In the case of nonnegative
potentials, they also proved an optimal upper estimate
$\frac{\lambda _{n}}{\lambda _{1}}\leq n^{p}$.

In this paper, we prove that if $q(x)\leq 0$ continuous\ and
single-barrier on the interval $[0,1],$ then $\frac{\lambda
_{n}}{\lambda _{m}}\geq \frac{n^{p}}{m^{p}}$ for $\lambda
_{n}>\lambda _{m}\geq -2q^{\ast },$ where $q^{\ast }=\inf \{ q( 0)
,q( 1) \} .$ Furthermore, we show that there exists $\ell_{0}\in (
0,1]$ such that for all $\ell \in(0,\ell_{0}],$ the associated
eigenvalues $ \lambda _{n}( \ell)$ $($of Problem $\eqref{1.1}
-\eqref{1.2} $ defined on $[ 0,\ell])$ satisfy $ \lambda _{1}(
\ell)>0$ and $\frac{\lambda _{n}( \ell)}{\lambda _{m}( \ell) }\geq
\frac{n^{p}}{m^{p}}$ for $n>m\geq 1$. The value $\ell _{0}$
satisfies the following estimate $0<\ell _{0}\leq
\sqrt[p]{\frac{-p}{3q^{*}}}$. As noted in Remark $\ref{rem1}$,our
approach used in this paper can be applied to the case of positive
potentials studied in $\cite{12}$ .

\section{Preliminaries And The  Main Statements }

In this section, Following $\cite{6}$, we introduce the modified
Pr\"{u}fer transformation.

Fix $p>1$.When $q \equiv 0$, equation $\eqref{1.1}$ becomes

\begin{eqnarray*}
-((y^{\prime }) ^{( p-1)})^{\prime }=(p-1)\lambda y^{(p-1)},x\in[
0,1]
\end{eqnarray*}
If $\lambda=1$ let $S_{p}$ be the solution of this equation
satisfying the initial conditions $S_{p}(0)=0,~~S'_{p}(0)=1.$ {A}.
Elbert $\cite{9} $, O. Dosly and P. Reh\'{a}k \cite{B}, proved that
$S_{p}(x)$ and $S'_{p}(x)$ are  $2\pi _{p}-$ periodic odd functions
on $\mathbb{R}$ (where $\pi_{p}=\frac{2\pi}{p\sin (\frac{\pi
}{p})}$). They are in fact $p$-analogues of $\sin $ and $\cos $
functions in the classical case. It is well known that $\pi
_{p}=\frac{2\pi }{p\sin ( \frac{\pi }{p}) }$ is the first zero of
$S_{p}$. The following lemma was proved by A. Elbert $\cite{9} $.

\begin{Lem}\label{lem1}$\cite{9} $ \begin{itemize}
                                       \item $\vert S_{p}(x)\vert ^{p}+\vert
                                       S_{p}^{\prime }(x)\vert ^{p}=1$ for any $x\in \R$.
                                       \item $\vert S_{p}^{\prime }(x)\vert
                                       ^{p-2}S_{p}^{"}(x)=-S_{p}^{( p-1) }( x) $ for any $x\in \R$.
                                   \end{itemize}

\end{Lem}
Let $T_{p}( x) =\frac{S_{p}(x)}{S_{p}^{\prime }(x)}.$ Then one can
easily derive $T_{p}^{\prime }( x) =1+\vert T_{p}( x) \vert ^{p}$
for any $x\in \R.$
\newline
We consider a generalized Pr\"{u}fer substitution on the
one-dimensional p-Laplacian $\eqref{1.1}$. With $\rho =\lambda
^{\frac{1}{p}}$ we have
\begin{equation}
y(x)=R(x)S_{p}(\varphi ( x) ),\ \ \ y^{\prime }(x)=\rho
R(x)S_{p}^{\prime }(\varphi ( x) )  \label{PR}
\end{equation}
\begin{eqnarray*}
\varphi ( 0,\rho ) =0,
\end{eqnarray*}%
where $R(x)>0$, and "prime" $({resp.~ "dot"})$\ denotes the
derivative with respect to $x$ $({resp. ~z}).$ Define further
\begin{eqnarray*}
\theta =\frac{\varphi }{\rho }.
\end{eqnarray*}
\begin{Lem}\label{lem2}$\cite{12,13}$
\begin{itemize}
                                       \item \begin{equation} \varphi ^{\prime }( x)
                                        =\rho - \frac{q( x)}{\rho^{p-1}} \vert S_{p}(\varphi ( x) )\vert ^{p} \label{teta},
                                        \end{equation}%
\item \begin{equation} \frac{R^{\prime}(x)}{R(x)}=\frac{q( x)}{\rho^{p-1}}S_{p}(\varphi ( x)
)^{p-1}S_{p}^{\prime }(\varphi (x)) \label{R.R}.
\end{equation}
 \end{itemize}
\end{Lem}
In the sequel we enunciate the main results of this paper.
\begin{Theo}\label{the1}
Let $q(x)\leq 0$ be monotone increasing in $[ 0,x_{0}] .$ Then $%
\dot{\theta}( x_{0},\rho ) \leq 0$ for\ $\rho \geq \sqrt[p]{%
-2q( 0) }.$ If there is a $\rho \geq \sqrt[p]{-2q( 0) } $ with
$\dot{\theta}( x_{0},z) =0,$ then $q\equiv0$ in $(0,x_{0}]$.
\end{Theo}
\begin{Theo}\label{the2}
For the system $ \eqref{1.1}- \eqref{1.2} ,$ if $q(x)\leq 0,$ continuous and single-barrier, then the $%
n^{th}$ and $m^{th}$ eigenvalues with $\lambda _{n}>\lambda _{m}\geq
-2q^{\ast }$ $($where $q^{\ast }=\inf \{ q( 0), q( 1) \} )$ satisfy
\begin{eqnarray*}
\frac{\lambda _{n}}{\lambda _{m}}\geq \frac{n^{p}}{m^{p}},
\end{eqnarray*}
and if for two different $m$ and $n$ equality holds, then $q\equiv 0
$ in $( 0,1).$
\end{Theo}
Consider the following boundary problem
\begin{equation}\label{1}
-((y^{\prime }) ^{( p-1)})^{\prime }=(p-1)(\lambda
-q(x))y^{(p-1)},x\in[ 0,\ell]
\end{equation}
\begin{equation}
y(0)=y(\ell )=0,  \label{T}
\end{equation}%
and let $\lambda _{n}( \ell ) ,$ $n\geq 1$ be the associated
eigenvalues.
\begin{Theo}\label{the3}
Assume that $q\leq 0,$ continuous and single-barrier in $[ 0,1] $.
Then, there exists $\ell _{0}\in ( 0,1] $ such that for all
$\ell\in(0,\ell_{0}]$ the associated eigenvalues $\lambda _{n}(
\ell)$ $($$of ~Problem ~ \eqref {1} -\eqref{T} ~defined ~on~
[0,\ell])$ satisfy $\lambda _{1}( \ell)>0$ and for the $n^{th}$ and
$m^{th}$ eigenvalues with $n>m\geq 1,$
\begin{eqnarray*}
\frac{\lambda _{n}( \ell)}{\lambda _{m}( \ell)}\geq
\frac{n^{p}}{m^{p}},
\end{eqnarray*}%
and if for two different $m$ and $n$ equality holds, then $q\equiv 0$ in $%
[ 0,\ell _{0}].$ \\ The value $\ell _{0}$ satisfies the following
estimate $0<\ell _{0}\leq \sqrt[p]{\frac{-p}{3q^{*}}}$, where
$q^{*}=\inf\{q(0) , q(1)\}.$
\end{Theo}

\begin{Rem}\label{rem1}
For $q(x) \geq 0$ and single-well we have $\varphi ^{\prime
}(x,\lambda )>0$ for $x\in[0,x_{0}]$ and $\lambda>0$. This implies
that the associated $\dot{\theta}(\lambda)$ can be written as entire
series in $\lambda>0$. Therefore, the ratios estimation
$\frac{\lambda_{n}}{\lambda_{m}}\leq (\frac{n}{m})^{p}$ will be
obtained for $n>m\geq1.$
\end{Rem}

Department of Mathematics. Faculty of Sciences Tunisia. University
EL-Manar Tunisia
\\ E-mail address: jamel.benamara@fsb.rnu.tn \\~~~~~~~~~~~\\
Department of Mathematics. Faculty of Sciences  Bizerte. University
of Carthage Tunisia. \\ E-mail address: hjihed@gmail.com

\begin{thebibliography}{99} \it{
\bibitem{9} \'{A}. Elbert: A half-linear second order differential equation.
Colloq. Math. Soc. J\'{a}nos Bolyai $30(1979),158-180.$
\bibitem{C} M. S. Ashbaugh and R. Benguria, "Best constant for the ratio of the
first two eigenvalues of onedimensional Schr¨odinger operators with
positive potentials,"Proc. Amer. Math. Soc., vol. $99, no. 3, pp.
598–599, 1987.$
\bibitem{1} M. S. Ashbaugh and R. D. Benguria, Optimal bounds for ratios of
eigenvalues of one dimensional Schr\"{o}dinger operators with
Dirichlet boundary conditions and positive potentials, Comm. Math.
Phys. $124$ $(1989)$ $419 - 424$
\bibitem{2} M. S. Ashbaugh and R. D. Benguria, Optimal lower bound for the
gap between the first two eigenvalues of one-dimensional Schrodinger
operators with symmetric single-well potentials, Proc. Amer. Math. Soc. $%
105, $ $419-424$ $(1989)$.
\bibitem{15} P. Binding and P. Dr\'{a}bek; Sturm-Liouville for the
p-laplacian. StudiaScientiarum Mathematicarum Hungaria $40$ $( 2003)
$ $375-396.$
\bibitem{4} Chung-Chuan Chen, Optimal lower estimates for eigenvalue ratios
of Schr\"{o}dinger operators and vibrating strings, electronically
published thesis $(2002)$.
\bibitem{12} Chao-Zhong Chen, et al.; Optimal Upper Bounds For The Eigenvalue Ratios of
One-Dimensional p-laplacian. Proc. Amer. Math. Soc. V $141,$N $3$,
$( 2013) $,$p$ $( 883-893) $.
\bibitem{6} J. B. Amara and J. Hedhly; A bound for ratios of eigenvalues of
Schr\"{o}dinger operators with single-barrier potentials. Submitted.
\bibitem{5} Mikl\'{o}s Horv\'{a}th and M\'{a}rton Kiss: A bound for ratios
of eigenvalues of Schr\"{o}dinger operators with single-well
potentials. Proc. Amer. Math. Soc. $134,$ $( 2005)$ $ 1425-1434.$
\bibitem{8} B. M. Levitan and I. S. Sargsjan, Sturm-Liouville and Dirac
operators (in Russian), Nauka, Moscow, $(1988).$
\bibitem{B} O. Do\~{s}ly and Reh\'{a}k, "Half-Linear Differential
Equations," Handbook of Differential Equations, North Holland
Mathematics Studies 202, Elsevier, Amsterdam, The Netherlands, $
2005.$
\bibitem{13} Gabriella Bogn\'{a}r and O. Do\~{s}ly, The ratio of eigenvalues of
the Dirichlet eigenvalue problem for equation with one-dimensional
p-Laplacian, Abstract and Applied Analysis $(2010)$
\bibitem{HL} Yu-Ling Huang and C. K. Law: Eigenvalue ratios for the regular
Sturm-Liouville system. Proc. Amer. Math. Soc. $124~(1996)~1427 -
1436.$
\bibitem{14} W. Walter, Sturm-Liouville theory for the radial $\triangle
_{p}-operator$, Math. Z., $227(1998),175-185.$}





\end{thebibliography}
\end{document}